 \documentclass[11pt]{article}
 \usepackage[top=1in,bottom=1in,left=1.5in,right=1.5in]{geometry}
  \usepackage{amsmath,amssymb}
  \usepackage{latexsym}% defines $\Box$ for LaTeX2e
\usepackage{graphicx} %SUS added
 \usepackage{subfigure} %SUS added
  \newcommand{\C}{\mathbb{C}}
  
  \newcommand{\N}{\mathbb{N}}
  
  \newcommand{\R}{\mathbb{R}}

  \newcommand{\bF}{\mathbf{F}}

  \newcommand{\m}{\mathbf{m}}

  \newcommand{\bt}{\mathbf{t}}

  \newcommand{\U}{\mathbf{U}}
  \newcommand{\uu}{\mathbf{u}}
  \newcommand{\vv}{\mathbf{v}}
  \newcommand{\V}{\mathbf{V}}
  \newcommand{\w}{\mathbf{w}}
  \newcommand{\W}{\mathbf{W}}
  \newcommand{\x}{\mathbf{x}}
  \newcommand{\X}{\mathbf{X}}
  \newcommand{\y}{\mathbf{y}}
  
  \newcommand{\z}{\mathbf{z}}
  \newcommand{\0}{\mathbf{0}}

  \newcommand{\bQ}{\mathbf{Q}}

  \newcommand{\cS}{\mathcal{S}}
  \newcommand{\cT}{\mathcal{T}}

  \newcommand{\cX}{\mathcal{X}}

  \newcommand{\lan}{\langle}
  \newcommand{\ran}{\rangle}
  \newcommand{\an}[1]{\lan#1\ran}
  
  \newcommand{\hs}{\hspace*{\parindent}}
  \newcommand{\proof}{\hs \textbf{Proof.\ }}

  \newcommand{\trans}{^\top}
  \newcommand{\qed}{\hspace*{\fill} $\Box$\\}

  \newcommand{\rS}{\mathrm{S}}

  \newtheorem{theo}{\bfseries \hs Theorem}%[section]
  
  \newtheorem{prop}[theo]{\bfseries \hs Proposition}
  
  \newtheorem{lemma}[theo]{\bfseries \hs Lemma}

  \numberwithin{equation}{section} % Automatically number equations within sections

 \renewcommand{\span}{\mathrm{span}}

\begin{document}

 \title{On best rank one approximation of tensors}

 \author{
 S. Friedland\footnotemark[1] \and
 V. Mehrmann\footnotemark[2]\and R. Pajarola\footnotemark[3] \and S.K. Suter\footnotemark[3]
 }
 \renewcommand{\thefootnote}{\fnsymbol{footnote}}
 \footnotetext[1]{
 Dept. of Mathematics, Statistics and Computer Science,
 Univ. of Illinois at Chicago, Chicago, Illinois 60607-7045,
 USA,
 \texttt{friedlan@uic.edu}.  This work was supported by NSF grant DMS-1216393.
 }
 \footnotetext[2]{
 Inst. f. Mathematik, MA4-5, TU Berlin, Str. des 17. Juni 136, D-10623 Berlin, FRG.
 \texttt{mehrmann@math.tu-berlin.de}. This work was supported by Deutsche Forschungsgemeinschaft (DFG) project  ME790/28-1.
 }
 \footnotetext[3]{Dept. of Informatics, Univ. of Z\"urich, Z\"urich, Switzerland \texttt{\{pajarola,susuter\}@ifi.uzh.ch}.
 Susanne Suter was partially supported by the Swiss National Science
 Foundation under Grant 200021-132521.}

 \renewcommand{\thefootnote}{\arabic{footnote}}
 \date{}
 \maketitle
\begin{abstract}
 In this paper we suggest a new algorithm for the computation of a best rank one approximation of tensors, called \emph{alternating singular value
 decomposition}.  This method is based on the computation of maximal singular values and the corresponding singular vectors of matrices.  We also
 introduce a modification for this method and the alternating least squares method, which ensures that alternating iterations will always converge to a
 semi-maximal point.  (A critical point in several vector variables is semi-maximal if it is maximal with respect to each vector variable, while other
 vector variables are kept fixed.)
% Finally, we introduce a new simple Newton-type method for speeding up the convergence of %alternating methods near the optimum.
 We present several numerical examples that illustrate the computational performance of the new method in comparison to the alternating least square
 method.
 \end{abstract}



 \renewcommand{\thefootnote}{\arabic{footnote}}

 \section{Introduction}\label{intro}
In this paper we consider the best rank one approximation to real
\emph{$d$-mode tensors} $\cT=[t_{i_1,\ldots,i_d}]\in\R^{m_1\times\ldots\times m_d}$,
 i. e.,  d-dimensional arrays with real entries.

As usual when studying tensors, it is necessary to introduce some notation.
Setting $[m]=\{1,\ldots,m\}$ for a positive integer $m$,
for two d-mode tensors $\cT,\cS \in\R^{m_1\times\ldots\times m_d}$ we denote by
\[
\an{\cT,\cS}:=\sum_{i_j\in[m_j], j\in [d]}t_{i_1,\ldots,i_d} s_{i_1,\ldots,i_d}
\]
the standard inner product of $\cT,\cS$, viewed as vectors in
$\R^{m_1\cdot m_2\cdot\ldots\cdot m_d}$. For an integer $p\leq d$, $r\in [p]$  and for $\x_{j_r}=[x_{1,j_r},\ldots,x_{m_{j_r},j_r}]\trans\in
\R^{m_{j_r}}$,
we use the standard mathematical notation
\[
\otimes_{j_r,r\in[p]}\x_{j_r}:=\x_{j_1}\otimes\ldots\otimes \x_{j_p}=[t_{i_1,\ldots,i_p}]\in\R^{m_{j_1}\times\ldots\times m_{j_p}}, \;
t_{i_1,\ldots,i_p}=x_{i_1,j_1} \ldots x_{i_p,j_p}.
\]
(See for example \cite[Chapter 5]{Fribk}.  In \cite{KolB09}  $\x\otimes\y$ is denoted as $\x\circ \y$ and is called vector outer product.)

 For a subset $P=\{j_1,\ldots,j_p\}\subseteq\emph{}[d]$ of cardinality $p=|P|$,
 consider a $p$-mode tensor $\cX=[x_{i_{j_1},\ldots,i_{j_p}}]\in\R^{m_{j_1}\times\ldots\times m_{j_p}}$, where $j_1<\ldots <j_p$.
 Then we have that $\cT\times \cX:=\sum_{i_{j_r}\in [m_{j_r}], r\in [p]} t_{i_1,\ldots,i_d} x_{i_{j_1},\ldots,i_{j_p}}$ is a $(d-p)$-mode tensor obtained
 by contraction on the indices $i_{j_1},\ldots,i_{j_p}$.  For example, if $\cT=[t_{i,j,k}]\in \R^{m\times n\times l}$ and $\x=[x_1,\ldots,x_{m}]\trans\in
 \R^m, \z=[z_1,\ldots,z_l]\trans \in \R^l$, then $\cT\times (\x\otimes\z):=\sum_{i\in[m],k\in[l]} t_{i,j,k}x_iz_k$, and it is viewed as a column vector in
 $\R^n$.
 Note that for  $\cT,\cS\in\R^{m_1\times\ldots\times m_d}$, we have $\an{\cT,\cS}=\cT\times\cS$.

For $\x\in \R^n$ we denote by $\|\x\|$ the Euclidian norm and
for $A\in \R^{m\times n}$ by $\|A\|=\max_{\|\x\|=1} \| A\x\|$ the associated operator norm. Then it is well-known, see e. g. \cite{GolV96}, that the best
rank one approximation of $A$ is given by $\sigma_1 \uu_1 \vv_1^T$, where
$\sigma_1=\| A \|$ is the largest singular value of $A$, and $\uu_1,\vv_1$ are the associated left and right singular vectors. Since the singular vectors
have Euclidian norm $1$, we have
that the spectral norm of the best rank one approximation is equal to $\sigma_1=\|A\|$.

To extend this property to tensors, let us for simplicity of  exposition restrict ourselves in this introduction to
the case of $3$-mode tensors $\cT\in\R^{m\times n\times l}$. Denote by $\rS^{m-1}:=\{\x\in \R^m, \;\|\x\|=1\}$ the unit sphere in $\R^m$, by $\rS(m,n,l)$
the set $\rS^{m-1}\times\rS^{n-1}\times\rS^{l-1}$, and introduce for $(\x,\y,\z)\in\rS(m,n,l)$ the function $f(\x,\y,\z):=\an{\cT,\x\otimes\y\otimes\z}$.
Then computing the best rank one approximation to $\cT$ is equivalent to finding
 \begin{equation}\label{bestrank1apr}
 \max_{(\x,\y,\z)\in\rS(m,n,l)} f(\x,\y,\z)=f(\x_\star,\y_\star,\z_\star).
 \end{equation}
The tensor version of the singular value relationship
%critical points of $f(\x,\y,\z)$ are then corresponding \emph{singular vectors} %of $\cA$ \cite{Lim05}
takes the form, see \cite{Lim05},
\begin{equation}\label{3tensingvect}
 \cT\times (\y\otimes\z)=\lambda\x,\; \cT\times (\x\otimes\z)=\lambda\y,\;\cT\times (\x\otimes\y)=\lambda\z,
 \end{equation}
where $\|\x\|=\|\y\|=\|\z\|=1$ and  $\lambda$ is a singular value of $\cT$.

Let us introduce for $p\in\{1,2\}$ the concept of a \emph{$p$-semi-maximum of $f$} restricted to $\rS(m,n,l)$.
For $p=1$, the $p$-semi-maximal points $\x_*,\y_*,\z_*$ of $f$ are the global maxima for the three functions
$f(\x,\y_*,\z_*)$, $f(\x_*,\y,\z_*)$, $f(\x_*,\y_*,\z)$ restricted to $\rS^{m-1}$, $\rS^{n-1}$, $\rS^{l-1}$, respectively. For $p=2$, the $p$-semi maximal
points are the pairs  $(\y_*,\z_*)$, $(\x_*,\z_*)$, $(\x_*,\y_*)$ that are global maxima of the functions $f(\x_*,\y,\z)$, $f(\x,\y_*,\z)$,
$f(\x,\y,\z_*)$
on $\rS^{n-1}\times\rS^{l-1}$, $\rS^{m-1}\times\rS^{l-1}$, $\rS^{m-1}\times\rS^{n-1}$, respectively.
We call $(\x_*,\y_*,\z_*)$ a semi-maximum if it is a $p$-semi-maximum for
$p=1$ or $p=2$, and it is clear how this concept of $p$-semi-maxima extends to arbitrary d-mode tensors with $p=1,2,\ldots,d-1$.
In the Appendix we discuss in detail $1$-local semi-maximal points of functions.

Many approaches for finding the maximum in \eqref{bestrank1apr} have been studied in the literature, see e. g. \cite{KolB09}. An important method, the
standard \emph{alternating least square} (ALS) method,  is an iterative procedure that
starts with $\x_0\in\rS^{m-1},\y_0\in\rS^{n-1},\z_0\in\rS^{l-1}$, where $f(\x_0,\y_0,\z_0)\ne 0$ and then defines the iterates $\x_i,\y_i,\z_i$ via
 \begin{equation}\label{ALSit}
 \x_{i}=\frac{\cT\times (\y_{i-1}\otimes \z_{i-1})}{\|\cT\times (\y_{i-1}\otimes \z_{i-1})\|},\
 \y_{i}=\frac{\cT\times (\x_{i}\otimes \z_{i-1})}{\|\cT\times (\x_{i}\otimes \z_{i-1})\|},\
 \z_{i}=\frac{\cT\times (\x_{i}\otimes \y_{i})}{\|\cT\times (\x_{i}\otimes \y_{i})\|},
 \end{equation}
 for $i=1,2,\ldots,\ $.

Note that  for all $i\in\N$ we have
\[
f(\x_{i-1},\y_{i-1},\z_{i-1})\le f(\x_{i},\y_{i-1},\z_{i-1})\le f(\x_{i},\y_{i},\z_{i-1}) \le f(\x_{i},\y_{i},\z_{i}),
\]
i. e., $f(\x_{i},\y_{i},\z_{i})$ is monotonically increasing and thus converges to a limit, since $f$ is bounded.  Typically, $(\x_{i},\y_{i},\z_{i})$
will
converge to a semi-maximum $(\x,\y,\z)$ that satisfies \eqref{3tensingvect}, however this is not clear in general \cite{KolB09}.
To overcome this deficiency of the ALS and related methods is one of the results of this paper.

We first discuss an alternative to the ALS algorithm for finding the maximum \eqref{bestrank1apr}, where each time
we fix only one variable and maximize on the other two.  Such a maximization is equivalent to finding the maximal singular value and the corresponding
left and right singular vectors of a suitable matrix, which is a well-established computational procedure, \cite{GolV96}. We call this method the
\emph{alternating singular value decomposition} (ASVD).
Next we introduce modifications of both ALS and ASVD, that are computationally more expensive, but for which it is guaranteed that they will always
converge to a semi-maximum of $f$.

Our numerical experimentation do not show clearly that ASVD is always better than ALS.  Since the standard algorithm for computing
the maximal singular value of a matrix is a truncated SVD algorithm \cite{GolV96}, and not ALS, we believe that ASVD is a very valid option in
finding best rank one approximations of tensors.
%The first iteration of the modified methods are three times more expensive then %the corresponding ALS and ASVD.  All other iterations are twice more
%expensive.

%Finally, as alternative to the Newton method suggested in \cite{ZhaG01}, we also suggest a %new Newton-type method  for best rank one approximation
%that is based on a simple observation mentioned in \cite{Fri11}, that any triple of singular %vectors $(\x,\y,\z)\in\rS(m,n,l)$
%corresponding to a nonzero singular value $\lambda$ of a tensor $\cT$ as in %\eqref{3tensingvect}, can be easily transferred to a fixed point of the map
%
%\begin{equation}\label{defmapF}
%\bF(\uu,\vv,\w):=(\cT\times(\vv\otimes\w),\cT\times(\uu\otimes\w),\cT\times(\uu\otimes\vv)).
%\end{equation}
%
%from $\R^m\times\R^n\times\R^l$ to itself.  Indeed, \eqref{3tensingvect} is equivalent to
%
%\begin{equation}\label{fixpointF}
%\bF(\uu,\vv,\w)-(\uu,\vv,\w)=\0,
%\end{equation}
%
%where $(\uu,\vv,\w)=\frac{1}{\lambda}(\x,\y,\z)$.  (This observation holds for any $d$-tensor %for $d>2$ but not for matrices!)

The content of the paper is as follows. In section~\ref{sec:basfacts} we recall some basic facts about tensors and best rank one approximations.
In section~\ref{sec:alt} we recall the ALS method and introduce the
ASVD procedure. The modification of these methods to guarantee convergence to a
semi-maximum is introduced in section~\ref{sec:mod} and
%the modification of the Newton iteration in section~\ref{sec:newton}.  T
the performance of the
new methods is illustrated in section~\ref{sec:num}.   In section~\ref{sec:conc} we state the conclusions of the paper.
In an Appendix we discuss the notion of local semi-maximality, give examples and discuss conditions for which ALS converges to a
local semi-maximal point.

\section{Basic facts on best rank one approximations of $d$-mode tensors}\label{sec:basfacts}
In this section we present further notation and recall some known results
about best rank one approximations.

For  a $d$-mode tensor $\cT=[t_{i_1,\ldots,i_d}]\in\R^{m_1\times\ldots\times m_d}$, denote by $\|\cT\|:=\sqrt{\an{\cT,\cT}}$ the Hilbert-Schmidt norm.
Denote by $\rS(\m)$ the $d$-product of the sub-spheres $\rS^{m_1-1}\times\ldots\times\rS^{m_d-1}$, let $(\x_1,\ldots,\x_d)\in\rS(\m)$  and associate with
$(\x_1,\ldots,\x_d)$ the $d$ one dimensional subspaces $\U_i=\span(\x_i)$, $i\in [d]$.
Note that
\[
\|\otimes_{i\in[d]}\x_i\|=\prod_{i\in[d]}\|\x_i\|=1.
\]
The projection $P_{\otimes_{i\in[d]}\U_i}(\cT)$  of $\cT$ onto the one dimensional subspace $\U:=\otimes_{i\in[d]}\U_i \subset \otimes_{i\in[d]}\R^{m_i}$,
 is given by
\begin{equation}\label{projform}
f(\x_1,\ldots,\x_d)\otimes_{i\in[d]}\x_i, \quad f(\x_1,\ldots,\x_d):=\an{\cT,\otimes_{i\in[d]}\x_i},\;(\x_1,\ldots,\x_d)\in\rS(\m).
\end{equation}
Denoting by $P_{(\otimes_{i\in[d]}\U_i)^\perp}(\cT)$ the orthogonal projection of $\cT$ onto the orthogonal complement of $\otimes_{i\in[d]}\U_i$,
the Pythagoras identity yields that
\begin{equation}\label{pythiden}
\|\cT\|^2=\|P_{\otimes_{i\in[d]\U_i}}(\cT)\|^2 + \|P_{(\otimes_{i\in[d]}\U_i)^\perp}(\cT)\|^2.
\end{equation}
With this notation, the best rank one approximation of $\cT$ from $\rS(\m)$ is given by
\[%begin{equation}\label{bestrank1ap}
\min_{(\x_1,\ldots,\x_d)\in\rS(\m)}\min_{a\in\R} \|\cT-a\otimes_{i\in[d]}\x_i\|.
\]%end{equation}
Observing that
 \[
 \min_{a\in\R} \|\cT-a\otimes_{i\in[d]}\x_i\|=
 \|\cT-P_{\otimes_{i\in[d]\U_i}}(\cT)\|=
 \|P_{(\otimes_{i\in[d]}\U_i)^\perp}(\cT)\|,
 \]
it follows that the best rank one approximation is obtained by the minimization of $\|P_{(\otimes_{i\in[d]}\U_i)^\perp}(\cT)\|$.
 In view of \eqref{pythiden} we deduce that best rank one approximation is obtained by  the maximization of $\|P_{\otimes_{i\in[d]\U_i}}(\cT)\|$ and
finally, using \eqref{projform}, it follows that the best rank one approximation is given by
\begin{equation}\label{specnorm}
\sigma_1(\cT):=\max_{(\x_1,\ldots,\x_d)\in\rS(\m)} f(\x_1,\ldots,\x_d).
\end{equation}
Following the matrix case, in \cite{HilL10} $\sigma_1(\cT)$ is called the \emph{spectral norm} and  it is also shown that the computation  of
$\sigma_1(\cT)$ in general is NP-hard for $d>2$.

We will make use of the following result of \cite{Lim05}, where  we present the proof for completeness.
\begin{lemma}\label{critptsf}  For $\cT\in\R^{m_1\times\ldots\times m_d}$, the
critical points of $f|_{\rS(\m)}$, defined in \eqref{projform}, satisfy the
equations
\begin{equation}\label{singvalvecd}
 \cT\times(\otimes_{j\in[d]\setminus\{i\}}\x_j)=\lambda \x_i\
  \textrm{ for all }i\in[d],\ (\x_1,\ldots,\x_d)\in\rS(\m), \end{equation}
for some real $\lambda$.
\end{lemma}
\proof We need to find the critical points of $\an{\cT,\otimes_{j\in[d]}\x_j}$ where $(\x_1,\ldots,\x_d)\in \rS(\m)$.
Using Lagrange multipliers we consider the auxiliary function
\[
g(\x_1,\ldots,\x_d):=\an{\cT,\otimes_{j\in[d]}\x_j}-\sum_{j\in[d]}\lambda_j \x_j\trans\x_j.
\]
The critical points of $g$ then satisfy
\[
\cT\times(\otimes_{j\in[d]\setminus\{i\}}\x_j)=\lambda_i \x_i, \quad i\in [d],
\]
and hence
$\an{\cT,\otimes_{j\in[d]}\x_j}=\lambda_i\x_i\trans\x_i=\lambda_i$ for all $i\in[d]$ which implies \eqref{singvalvecd}.
\qed

Observe next that $(\x_1,\ldots,\x_d)$ satisfies  \eqref{singvalvecd} iff the vectors $(\pm\x_1,\ldots,\pm\x_d)$ satisfy \eqref{singvalvecd}.
 In particular, we could choose the signs in $(\pm\x_1,\ldots,\pm\x_d)$ such that each corresponding $\lambda$ is nonnegative and then these
$\lambda$ can be interpreted  as the singular values of $\cT$.  The maximal singular value of $\cT$ is denoted by $\sigma_1(\cT)$ and is given
 by \eqref{specnorm}.  Note that to each nonnegative singular value there are at least $2^{d-1}$ singular vectors of the form $(\pm\x_1,\ldots,\pm\x_d)$.
 So it is more natural to view the singular vectors as one dimensional subspaces $\U_i=\span(\x_i)$, $i\in[d]$.

 As observed in \cite{Fri11} for tensors, i. e., for $d>2$, there is a one-to-one correspondence between the singular vectors corresponding to positive
 singular values of $\cT$ and the fixed points of an induced multilinear map of degree $d-1$.
 \begin{lemma}\label{lemfixpointsF}  Let $d>2$  and assume that $\cT\in\R^{m_1\times\ldots\times m_d}$.  Associate with $\cT$ the map $\bF$ from
 $\R(\m):=\R^{m_1}\times\ldots\times\R^{m_d}$ to itself, where
\[%begin{equation}\label{defFd}
 \bF:=(F_1,\ldots,F_d), \;F_i(\uu_1,\ldots,\uu_d):=\cT\times(\otimes_{j\in[d]\setminus\{i\}} \uu_j), \;i\in[d].
\]%end{equation}
Then there is a one-to-one correspondence between the critical points of $f|_{\rS(\m)}$ corresponding to positive singular values $\lambda$ and the
nonzero
fixed points of
\begin{equation}\label{deffixptsF}
 \bF(\uu)=\uu.
\end{equation}
Namely, each $(\x_1,\ldots,\x_d)\in\rS(\m)$ satisfying \eqref{singvalvecd} with $\lambda>0$ induces a fixed point of $\bF$
of the form
\[
(\uu_1,\ldots,\uu_d)=\lambda^{\frac{-1}{d-2}}(\x_1,\ldots,\x_d).
\]
Conversely, any nonzero fixed point satisfying \eqref{deffixptsF}
induces a $d$-set of singular vectors $(\x_1,\ldots,\x_d)=\frac{1}{\|\uu_1\|}(\uu_1,\ldots,\uu_d)\in\rS(\m)$ corresponding to
$\lambda=\|\uu_1\|^{-(d-2)}$.
In particular, the spectral norm $\sigma_1(\cT)$ corresponds to a nonzero fixed point of $\bF$ closest to the origin.
\end{lemma}
\proof
Assume that  \eqref{singvalvecd} holds for $\lambda>0$.  Dividing both sides of \eqref{singvalvecd} by $\lambda^{\frac{d-1}{d-2}}$ we obtain that
$(\uu_1,\ldots,\uu_d)=\lambda^{\frac{-1}{d-2}}(\x_1,\ldots,\x_d)$ is a nonzero fixed point of $\bF$.

For the converse, assume that $(\uu_1,\ldots,\uu_d)$ is a nonzero fixed point of $\bF$.  Clearly $\uu_i\trans\uu_i=\an{\cT,\times_{j\in[d]}\uu_j}$ for
$i\in[d]$.  Hence, $\|\uu_1\|=\ldots=\|\uu_d\|>0$ and  $(\x_1,\ldots,\x_d)=\frac{1}{\|\uu_1\|}(\uu_1,\ldots,\uu_d)\in\rS(\m)$ satisfies
\eqref{singvalvecd}
with $\lambda=\|\uu_1\|^{-(d-2)}$.

The largest positive singular value of $\cT$ corresponds to the nonzero fixed point $(\uu_1,\ldots,\uu_d)$, where $\sum_{i\in[d]}\|\uu_i\|^2=d\|\uu_1\|^2$
is the smallest.
\qed
We also have that the trivial fixed point is isolated.
\begin{prop}\label{zerofixptF} The origin $\0\in\R(\m)$ is an isolated simple fixed point of $\bF$.
 \end{prop}
\proof  A fixed point of $\bF$ satisfies
\begin{equation}\label{deffixptsFm}
\uu-\bF(\uu)=\0
\end{equation}
and clearly, $\uu=\0$ satisfies this system.  Observe that the Jacobian matrix $D(\uu-\bF(\uu))(\0)$ is the identity matrix.
 Hence the implicit function theorem yields that $\0$ is a simple isolated solution of \eqref{deffixptsF}.
\qed

In view of Lemma~\ref{lemfixpointsF} and Proposition~\ref{deffixptsFm}, to compute the best rank one tensor approximation,
we will introduce an iterative procedure that converges to the fixed point closest to the origin.

In \cite{FO12} the following results are established.  First, for a generic $\cT\in\R^{m_1\times\ldots\times m_d}$
the best rank one approximation of $\cT$ is unique. Second, a complex generic $\cT\in\C^{m_1\times\ldots\times m_d}$
has a finite number $\nu(m_1,\ldots,m_d)$ of singular value tuples and the corresponding ``singular complex values" $\lambda$.
We now consider the ``cube" case where $m_1=\ldots=m_d=m$.  Then
$\nu(m,\ldots,m)$ is different from the number of complex eigenvalues computed in \cite{CS}.
Finally,  for a generic symmetric tensor $\cT\in \R^{m\times\ldots\times m}$, the best rank one approximation is unique and symmetric.
(The fact that the best rank one approximation of a symmetric tensor can be chosen symmetric is proved in \cite{Fri11}.)
\section{The ALS and the ASVD method}\label{sec:alt}
In this section we briefly recall the alternating least squares (ALS) method and suggest an analogous alternating singular value decomposition (ASVD)
method.

Consider $\cT\in\R^{m_1\times\ldots\times m_d}\setminus\{0\}$ and
choose an initial point $(\x_{1,0},\ldots,\x_{d,0})\in\rS(\m)$ such that $f(\x_{1,0},\ldots,\x_{d,0})\ne 0$. This can be done in different ways. One
possibility is to choose  $(\x_{1,0},\ldots,\x_{d,0})\in\rS(\m)$ at random.  This will ensure that with probability one we have
$f(\x_{1,0},\ldots,\x_{d,0})\ne 0$.
Another, more expensive way to obtain such an initial point $(\x_{1,0},\ldots,\x_{d,0})$ is to use the higher order singular value decomposition (HOSVD)
\cite{DeLDV00}.
 To choose $\x_{i,0}$ view $\cT$ as an $m_i\times \frac{m_1\times\ldots\times m_d}{m_i}$ matrix $A_i$, by unfolding in direction $i$.  Then $\x_i$ is the
 left
 singular vector corresponding to $\sigma_1(A_i)$ for $i\in[d]$.  The use of
 the HOSVD is expensive, but may result in a better choice of the initial point.

Given $(\x_{1,p},\ldots,\x_{d,p})\in\rS(\m)$, for an integer $p\ge 0$ the points $\x_{i,p+1}\in \rS^{m_i-1}$ are then computed recursively via
\begin{equation}\label{als}
\x_{i,p+1}=\frac{1}{\|\cT
\times(\otimes_{j=1}^{i-1}\x_{j,p+1}\otimes(\otimes_{j=i+1}^{d}\x_{j,p}) )\|}
(\cT\times((\otimes_{j=1}^{i-1}\x_{j,p+1})\otimes(\otimes_{j=i+1}^{d}\x_{j,p} ))),
\end{equation}
for $i\in[d]$. Each iterate of (\ref{als}) is the solution of an optimization problem which is obtained by setting
the gradient of a simple Lagrangian to $0$. % and $p\ge 0$.
Therefore, clearly, we have the inequality
\[
%\begin{equation}\label{firstinq}
f(\x_{1,p+1},\ldots,\x_{i-1,p+1}, \x_{i,p},\ldots,\x_{d,p})\le f(\x_{1,p+1},\ldots,\x_{i,p+1}, \x_{i+1,p},\ldots,\x_{d,p}),
%\end{equation}
\]
for $i\in[d]$ and $p\ge 0$, %.  Since
%
%\begin{eqnarray*}
% f(\x_{1,p+1},\x_{2,p+1}, %\x_{1,p},\ldots,\x_{d,p})&:=&f(\x_{1,p},\ldots,\x_{d,p}),\\
% f(\x_{1,p+1},\ldots,\x_{d,p+1}, %\x_{d+1,p},\ldots,\x_{d,p})&:=&f(\x_{1,p+1},\ldots,\x_{d,p+1}),
% \end{eqnarray*}
%
and the sequence $f(\x_{1,p},\ldots,\x_{d,p}), p=0,1,\ldots$ is a nondecreasing sequence bounded by $\sigma_1(\cT)$, and hence it converges.

Recall that a point $(\x_{1,*},\ldots,\x_{d,*})\in\rS(\m)$ is called a $1$-semi maximum, if $\x_{i,*}$ is a maximum for the function $f(\x_{1,*},\ldots,
\x_{i-1,*},\x_i,\x_{i+1,*},\ldots,\x_{d,*})$ restricted to $\rS^{m_i-1}$ for each $i\in[d]$. Thus, clearly any $1$-semi maximal point of $f$ is a critical
point of $f$.
In many cases the sequence $(\x_{1,p},\ldots,\x_{d,p}), p=0,1,\ldots$ does converge to a $1$-semi maximal point of $f$,
however, this is not always guaranteed \cite{KolB09}.
%however, there exist examples in which this not the case \cite{KolB09}.

An alternative to the ALS method is the alternating singular value decomposition (ASVD). To introduce this method, denote for $A\in\R^{m\times \ell}$ by
$\uu(A)\in\rS^{m-1},\vv(A)\in \rS^{\ell-1}$ the left and the right singular vectors of $A$
 corresponding to the maximal singular value $\sigma_1(A)$, i. e.,
\[
\uu(A)\trans A=\sigma_1(A)\vv(A)\trans, A\vv(A)=\sigma_1(A)\uu(A).
\]
Since for $d=2$ the singular value decomposition directly gives the best rank one approximation, we only consider the case $d\geq 3$.
Let $\cT=[t_{i_1,\ldots,i_d}]\in\R^{m_1\times\ldots\times m_d}$ and $X:=(\x_1,\ldots,\x_d) \in \rS(\m)$ be such  that $f(\x_1,\ldots,\x_d)\ne 0$.
Fix an index pair $(i,j)$ with $1\le i<j\le d$ and denote  by $\cX_{i,j}$ the $d-2$ tensor $\otimes_{k\in[d]\setminus\{i,j\}}\x_k$.  Then $\cT\times
\cX_{i,j}$ is an $m_i\times m_j$ matrix.

The basic step in the ASVD method is the substitution
\begin{equation}\label{replacestepij}
 (\x_i,\x_j)\mapsto (\uu(\cT\times\cX_{i,j}),\vv(\cT\times \cX_{i,j})).
\end{equation}
For example, if $d=3$ then the ASVD method is given by repeating iteratively the  substitution (\ref{replacestepij}) in the order
\[%begin{equation}\label{d=3pairord}
 (2,3),\quad (1,3),\quad (1,2).
\]%end{equation}
For $d>3$, one goes consecutively through all $d \choose 2$ pairs in an ``evenly distributed way".  For example, if $d=4$ then one could choose the order
\[
% \begin{equation}\label{d=4pairord}
 (1,2),\quad (3,4),\quad (1,3),\quad (2,4), \quad (1,4), \quad (2,3).
\]
% \end{equation}

Observe that the substitution \eqref{replacestepij} gives $\sigma_1(\cT\times \cX_{i,j})$.
Note that the ALS method for the bilinear form $g(\x,\y)=\x\trans (\cT\times \cX_{i,j})\y$ could increase the value of $g$ at most to its
maximum, which is $\sigma_1(\cT\times \cX_{i,j})$.  Hence we have the following proposition.
\begin{prop}\label{ASVDcomp}  Let  $\cT\in\R^{m_1\times\ldots\times m_d}\setminus\{0\}$ and  assume that
$(\x_{1},\ldots,\x_{d})\in\rS(\m)$. Fix $1\le i < j\le d$ and consider the following three maximization problems.
First, fix all variables except the variable $\x_p$ and denote the maximum of $f(\x_1,\ldots,\x_d)$ over $\x_p\in \rS^{m_p-1}$ by
$a_p$. Then find $a_i,a_j$.  Next fix all the variables except $\x_i,\x_j$  and find the maximum of $f(\x_1,\ldots,\x_d)$ over
$(\x_i,\x_j)\in \rS^{m_i-1}\times \rS^{m_j-1}$,  which is denoted  by $b_{i,j}$.  Then $b_{i,j}\ge \max(a_i,a_j)$.
In particular one step in the ASVD increases the value of $f$ as least as much as a corresponding step of ALS.
\end{prop}
The procedure to compute the largest singular value of a large scale
matrix is based on the Lanczos algorithm \cite{GolV96} implemented
in the partial singular value decomposition. Despite the fact
that this procedure is very efficient,
% Note that for a given matrix $A$  we do only one iteration of the form %\eqref{replacestepij} to find $\sigma_1(A)$, using partial SVD.  This is
%considered the most efficient way to compute $\sigma_1(A)$.
for tensors each iteration of ALS is still much cheaper to perform than one iteration of \eqref{replacestepij}.  However, it is not really necessary
to iterate the partial SVD algorithm to full convergence of the largest singular
value.  Since the Lanczos algorithm converges rapidly \cite{GolV96},
a few steps (giving only a rough approximation) may be enough to get an improvement in the outer iteration. In this case, the ASVD method may even  be
faster than the ALS method, however, a complete analysis of such an inner-outer iteration is an open problem.  As in the ALS method, it may happen that a
step of the ASVD will not decrease the value of the function $f$, but in many cases the algorithm will converge to a semi-maximum of $f$.  However, as in
the case of the ALS method, we do not have a complete understanding when this will happen.
For this reason, in the next section we suggest a modification of both ALS and ASVD method, that will guarantee convergence.

\section{Modified ALS and ASVD}\label{sec:mod}
The aim of this section is to introduce \emph{modified} ALS and ASVD methods, abbreviated here as MALS and MASVD.  These modified algorithms ensure that
every accumulation point of these algorithms is a semi-maximal point of $f|_{\rS(\m)}$.
%In particular, if one of these accumulation points is an isolated semi-maximum of $f|_{\rS(\m)}$ then the iteration will converge to a %semi-maximum.
For simplicity of the exposition we describe the concept for the case $d=3$, i. e., we assume that we have a tensor $\cT\in \R^{m\times n\times l}$.

 We first discuss the MALS.
For given $(\x,\y,\z)\in\rS(m,n,l)$ with $f(\x,\y,\z)\ne 0$, the procedure
 requires to compute the three  values
\begin{eqnarray*}
 f_{1}(\x,\y,\z)&:=&f \left (\frac{\cT\times (\y\otimes \z)}{\|\cT\times (\y\otimes \z)\|},\y,\z\right ),\;\\
 f_2(\x,\y,\z)&:=&f\left (\x,\frac{\cT\times (\x\otimes \z)}{\|\cT\times (\x\otimes \z)\|},\z\right ),\\
 f_3(\x,\y,\z)&:=&f\left (\x,\y,\frac{\cT\times (\x\otimes \y)}{\|\cT\times (\x\otimes \y)\|} \right),
\end{eqnarray*}
and to choose the maximum value. This needs  $3$ evaluations of $f$.

The modified ALS procedure then is as follows.
Let $(\x_0,\y_0,\z_0)\in\rS(m,n,l)$ and $f(\x_0,\y_0,\z_0)\ne 0$.
 Consider the maximum value of $f_i(\x_0,\y_0,\z_0)$ for $i=1,2,3$.  Assume for example that this value is achieved for $i=2$ and  let
 $\y_1:=\frac{\cT\times (\x_0\otimes \z_0)}{\|\cT\times (\x_0\otimes \z_0)\|}$.  Then we replace the point $(\x_0,\y_0,\z_0)$ with the new point
 $(\x_0,\y_1,\z_0)$ and
consider the maximum value of $f_i(\x_0,\y_1,\z_0)$ for $i=1,2,3$.
This needs only two $f$ evaluations, since $f_2(\x_0,\y_0,\z_0)=f_2(\x_0,\y_1,\z_0)$. Suppose that this maximum is achieved for $i=1$.  We then replace
the
point in the triple $(\x_0,\y_1,\z_0)$ with $(\x_1,\y_1,\z_0)$
where $\x_1=\frac{\cT\times (\y_1\otimes \z_0)}{\|\cT\times (\y_1\otimes \z_0)\|}$ and then as the last step we optimize the missing mode, which is in
this
example $i=3$. In case that the convergence criterion is not
yet satisfied, we continue iteratively in the same  manner.  The cost of
this algorithm is about twice as much as that of ALS.

For the modified ASVD we have a similar algorithm.  For $(\x,\y,\z)\in\rS(m,n,l)$, $f(\x,\y,\z)\ne 0$, let
 \begin{eqnarray*}
 g_1(\x,\y,\z)&:=& f(\x,\uu(\cT\times\x),\vv(\cT\times\x)),\\
  g_2(\x,\y,\z)&:=& f(\uu(\cT\times\y),\y,\vv(\cT\times\y)),\\
 g_3(\x,\y,\z)&:=& f(\uu(\cT\times\z),\vv(\cT\times\z),\z),
 \end{eqnarray*}
which requires three evaluations of $f$.
Let $(\x_0,\y_0,\z_0)\in\rS(m,n,l)$ and $f(\x_0,\y_0,\z_0)\ne 0$ and
consider the maximal value of $g_i(\x_0,\y_0,\z_0)$ for $i=1,2,3$.  Assume for example that this value is achieved for $i=2$.  Let $\x_1:=\uu(\cT\times
\y_0), \z_1:=\vv(\cT\times \y_0)$.  Then we replace the point $(\x_0,\y_0,\z_0)$ with the new point $(\x_1,\y_0,\z_1)$ and determine the maximal value of
$g_i(\x_1,\y_0,\z_1)$ for $i=1,2,3$.
Suppose that this maximum is achieved for $i=1$.  We then replace the point in the triple $(\x_1,\y_0,\z_1)$ with $(\x_1,\y_1,\z_2)$
where $\y_1=\uu(\cT\times\x_1), \z_2=\vv(\cT\times\x_1)$ and if the convergence criterion is not met then we continue in the same manner.  This algorithm
is about twice as expensive as the ASVD method. For $d=3$, we then have the following theorem.
\begin{theo}\label{convmodmeth}  Let $\cT\in \R^{m\times n\times l}$ be a given
tensor and consider the sequence
\begin{equation}\label{iterseq}
 (\x_i,\y_i,\z_i)\in \rS(m,n,l) \textrm{ for } i=0,1,\ldots,
\end{equation}
generated either by MALS or MASVD, where $f(\x_0,\y_0,\z_0)\ne 0$.  If $(\x_{*},\y_*,\z_*)\in \rS(m,n,l)$ is an accumulation point of this sequence, then
$(\x_{*},\y_*,\z_*)\in \rS(m,n,l)$ is a $1$-semi maximum if \eqref{iterseq} is given by MALS and a $2$-semi maximum if \eqref{iterseq} is given by MASVD.
 \end{theo}
\proof  Let $(\x_{*},\y_*,\z_*)\in \rS(m,n,l)$ be an accumulation point of the sequence \eqref{iterseq}, i.e., there exists a subsequence
$1\le n_1<n_2<n_3<\ldots$ such that

\noindent
$\lim_{j\to\infty} (\x_{n_j},\y_{n_j},\z_{n_j})=(\x_*,\y_*,\z_*)$.
Since the sequence $f(\x_i,\y_i,\z_i)$ is nondecreasing, we deduce that
$\lim_{i\to\infty} f(\x_i,\y_i,\z_i)=f(\x_*,\y_*,\z_*)>0$.
By the definition of $f_i(\x_*,\y_*,\z_*)$ it follows that
\begin{equation}\label{fiineqxyz}
 \min\{f_j(\x_*,\y_*,\z_*),\; j=1,2,3\}\ge f(\x_*,\y_*,\z_*).
\end{equation}
Assume first that the sequence \eqref{iterseq} is obtained by either ALS and MALS.
We will point out exactly, where we need the assumption that \eqref{iterseq} is obtained by MALS to deduce that  $(\x_{*},\y_*,\z_*)\in \rS(m,n,l)$ is a
$1$-semi maximum.

Consider first the ALS sequence given as in \eqref{ALSit}. Then
\begin{eqnarray}\notag
 &&f(\x_i,\y_{i-1},\z_{i-1})= f_1(\x_{i-1},\y_{i-1},\z_{i-1}) \\
 &&\le f(\x_i,\y_{i},\z_{i-1})= f_2(\x_i,\y_{i-1},\z_{i-1})\notag \\
 &&\le f(\x_i,\y_i,\z_i)=f_3(\x_{i},\y_{i},\z_{i-1}).\label{behASL}
\end{eqnarray}
For any $\varepsilon>0$, since $f_1(\x,\y,\z)$ is a continuous function on $\rS(m,n,l)$, it follows that for a sufficiently large  integer $j$  that
$f_1(\x_{n_j},\y_{n_j},\z_{n_j})>f_1(\x_*,\y_*,\z_*)-\varepsilon$.  Hence
\begin{equation}\label{fepsin}
 f(\x_*,\y_*,\z_*)\ge f(\x_{n_j+1},\y_{n_j+1},\y_{n_j+1})\ge f_1(\x_{n_j+1},\y_{n_j},\z_{n_j})\ge f_1(\x_*,\y_*,\z_*)-\varepsilon.
\end{equation}
Since $\varepsilon>0$ can be chosen arbitrarily small, we can combine inequality (\ref{fepsin})
with \eqref{fiineqxyz} to deduce that $f_1(\x_*,\y_*,\z_*)=f(\x_*,\y_*,\z_*)$.
We can also derive the equality  $f_3(\x_*,\y_*,\z_*)=f(\x_*,\y_*,\z_*)$ as follows.  Clearly,
\[
f(\x_{n_j},\y_{n_j},\z_{n_j-1})\le f_3(\x_{n_j},\y_{n_j},\z_{n_j-1})=f(\x_{n_j},\y_{n_j},\z_{n_j})\le f(\x_{n_{j+1}},\y_{n_{j+1}},\z_{n_{j+1}})
\]
Using the same arguments as for $f_1$ we deduce the equality $f_3(\x_*,\y_*,\z_*)= f(\x_*,\y_*,\z_*)$.
 However, there is no way to deduce equality in the inequality $f_2(\x_*,\y_*,\z_*)\ge f(\x_*,\y_*,\z_*)$ for the ALS method, %  as described above,
since $f_2(\x_{i},\y,\z_{i})=f(\x_i,\uu_i,\z_i)$  and $\uu_i$ is not equal to $\y_{i}$ or $\y_{i+1}$.

We now consider the case of MALS.  We always have the inequalities
$f_j(\x_i,\y_i,\z_i)\le f(\x_{i+1},\y_{i+1}, \z_{i+1})$ for each $j=1,2,3$ and $i\in\N$.  Then the same arguments as before imply in a straightforward way
that
we have equalities in \eqref{fiineqxyz}.  Hence $(\x_*,\y_*,\z_*)$ is  a $1$-semi maximum.

 Similar arguments show that if the sequence \eqref{iterseq} is obtained by MASVD then $g_k(\x_*,\y_*,\z_*)=f(\x_*,\y_*,\z_*)$ for $k\in[3]$.
 Hence $(\x_*,\y_*,\z_*)$ is a $2$-semi maximum.
\qed
It is easy to accelerate the convergence of the MALS and MASVD algorithm in the neighborhood of a semi-maximum via  Newton's method, see e.g.
\cite{ZhaG01}.

Despite the fact Theorem~\ref{convmodmeth} shows convergence to $1$- or $2$-semi-maximal points, the monotone convergence can not be employed to show convergence to a critical point and the following questions remain {\bf open}.  Suppose that the assumptions of Theorem~\ref{convmodmeth} hold.  Assume further, that one accumulation
point $(\x_{*},\y_*,\z_*)$ of the sequence \eqref{iterseq} is an isolated critical point of $f|_{\rS(m,n,l)}$.  Is it true that for the MALS method and a generic starting value the
sequence \eqref{iterseq} converges to $(\x_{*},\y_*,\z_*)$, where we identify $-\x_*,-\y_*,-\z_*$ with $\x_*,\y_*,\z_*$ respectively?  Is the same claim
true  for the MASVD method assuming the additional condition
 %
%\begin{equation}\label{nondegcondMASVD}
\[
 \sigma_1(\cT\times \x_*)>\sigma_2(\cT\times \x_*),\;\sigma_1(\cT\times \y_*)>\sigma_2(\cT\times \y_*),\;
 \sigma_1(\cT\times \z_*)>\sigma_2(\cT\times \z_*)?
\]
%\end{equation}
%
In the Appendix we show that for specific initial values convergence may not happen towards the unique isolated critical point, but towards other semi-maximal points. Our numerical results with random starting values however, seem to confirm
the convergence to the unique critical  point.

\section{Numerical results}\label{sec:num}
%\subsection{Simulation Setup}
We have implemented a C++ library supporting the rank one tensor decomposition using vmmlib~\cite{vmmlib}, LAPACK and BLAS in order to test the
performance
of the different best rank one approximation algorithms. The performance was measured via the actual CPU-time (seconds) needed to
compute the approximate best rank one decomposition,  by the number of optimization calls needed, and
whether a stationary point was found.
%whether a stationary point or a global maxima is found.
All performance tests have been carried out on a 2.8 GHz Quad-Core Intel Xeon Macintosh computer with 16GB RAM.

The performance results are discussed for synthetic and real data sets of third-order tensors. In particular, we worked with three different data sets:
(1)
a real computer tomography (CT) data set (the so-called MELANIX data set of OsiriX), (2) a symmetric random data set, where all indices are symmetric, and (3) a
random data set. The CT data set has a 16bit, the random data set an 8bit value range.
All our third-order tensor data sets are initially of size $512 \times 512 \times 512$, which we gradually reduced by a factor of $2$, with the smallest
data sets being of size $4 \times 4 \times 4$. The synthetic random data sets were generated for every resolution and in every run; the real data set was
averaged (subsampled) for every coarser resolution.

Our simulation results are averaged over different decomposition runs of the various algorithms. In each decomposition run, we changed the initial guess,
i.e., we generated new random start vectors. We always initialized the algorithms by random start vectors, since this is cheaper than the initialization
via HOSVD. Additionally, we generated for each decomposition run new random data sets. The presented timings  are averages over 10 different runs of the
algorithms.

All the best rank one approximation algorithms are alternating algorithms, and based on the same convergence criterion, where convergence is
achieved if one of the two following conditions: $iterations > 10$; $fitchange <  0.0001$ is met. The number of optimization calls within one iteration
is
fixed for the ALS and ASVD method. During one iteration, the ALS optimizes every mode once, while the ASVD optimizes every mode twice. The number of
optimization calls can vary widely during each iteration of the modified algorithms. This results from the fact that multiple optimizations are performed
in parallel, while only the best one is kept and the others are rejected.
The partial SVD is implemented by applying a symmetric eigenvalue decomposition (LAPACK DSYEVX) to the product $A A^T$ (BLAS DGEMM) as suggested by the
ARPACK package.

%SUS: started to edit from here

\begin{figure}[htbp]
\begin{center}
     \subfigure[CPU time (s) for medium sized 3-mode tensor samples]{\label{fig:times_st3}
     \includegraphics[width=0.9\columnwidth]{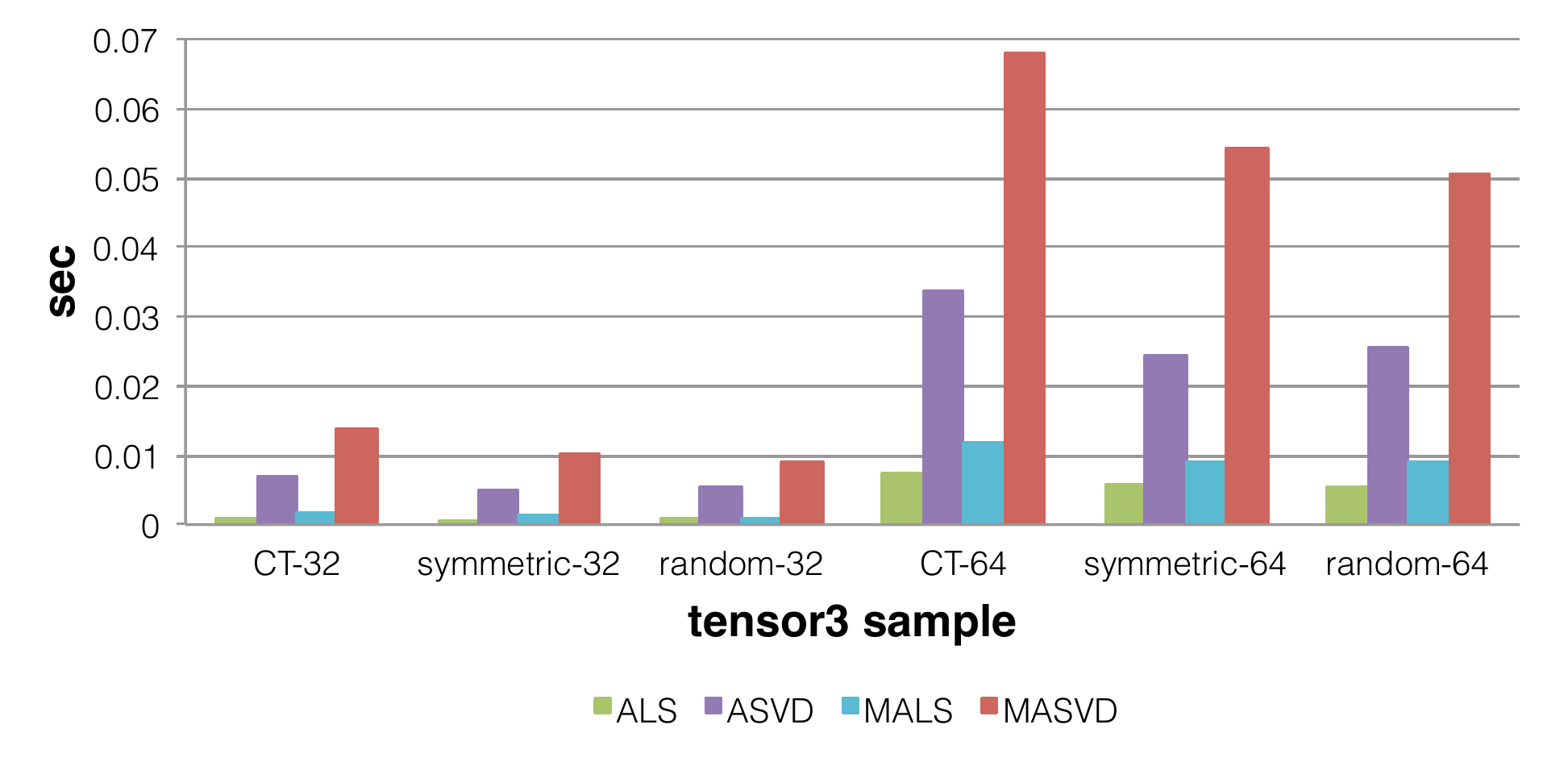}} \\
     \subfigure[CPU time (s) for larger 3-mode tensor samples]{\label{fig:times_lt3}
     \includegraphics[width=0.9\columnwidth]{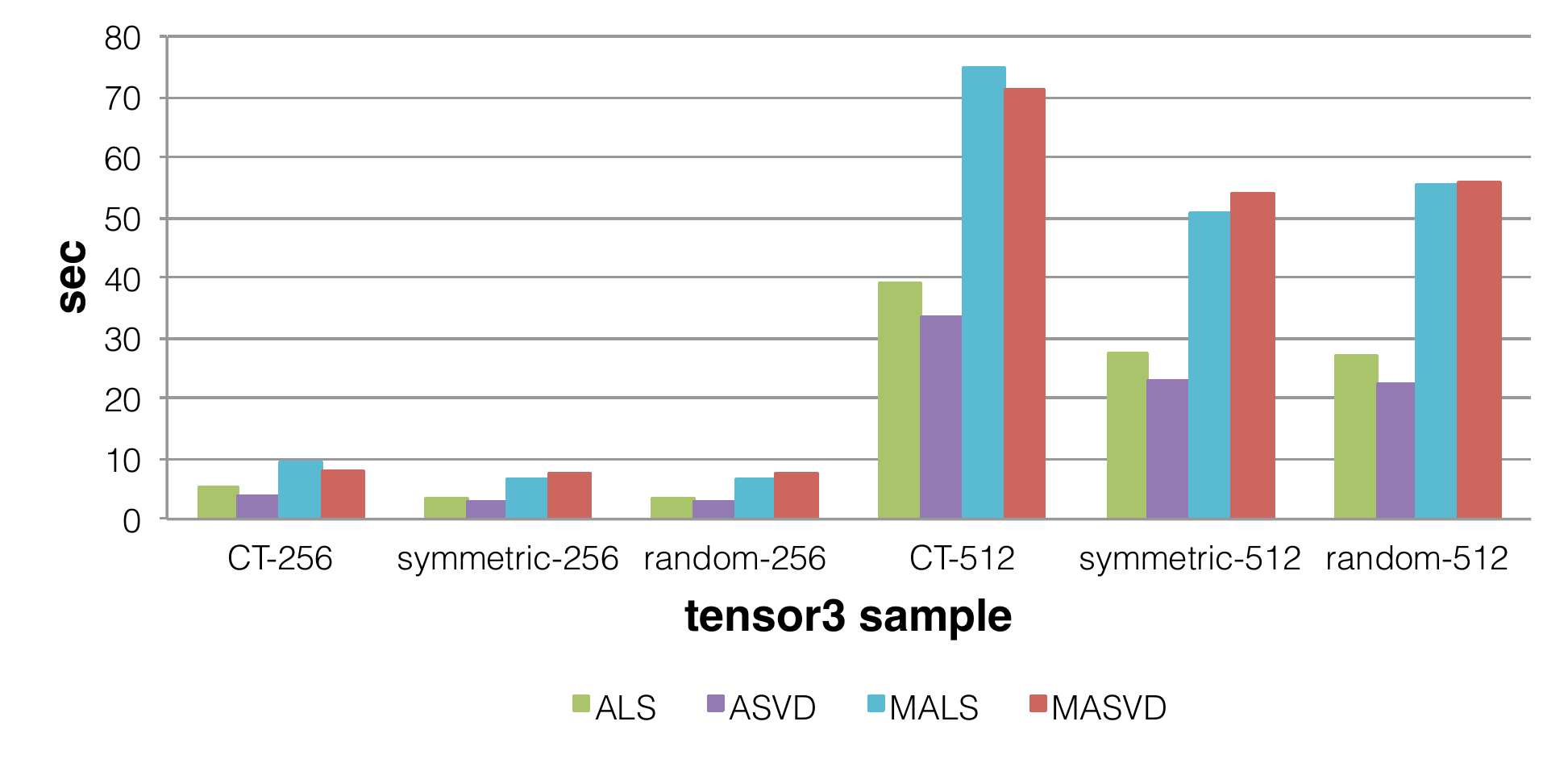}} \\
\caption{Average CPU times for best rank one approximations
per algorithm and per data set taken over 10 different initial random
guesses.}
\label{fig:times}
\end{center}
\end{figure}

With respect to the total decomposition times for different sized third-order tensors (tensor3s), we observed that for tensor3s smaller than $64^3$, the total decomposition time was below one second.
That represents a time range, where we do not need to optimize further.
On the contrary, the larger the tensor3s gets, the more critical the differences in the decomposition times are.
As shown in Figure~\ref{fig:times}, the modified versions of the algorithms took about twice as much CPU-time as the standard versions.
For the large data sets, the ALS and ASVD take generally less time than the MALS and MASVD. The ASVD was fastest
for large data sets, but compared to (M)ALS slow for small data sets. For larger data sets, the timings of the basic and modified algorithm versions came
closer to each other.

Furthermore, we compared the number of optimization calls needed for the algorithms of ALS, ASVD, MALS, and MASVD, recalling that for the ALS and the MALS, one
mode is optimized per optimization call, while for ASVD and MASVD, two modes are optimized per optimization call.
Figure 2 demonstrates the relationships of the four algorithms (color encoded) on three different data sets (marker encoded) and the different data set sizes (hue encoded). As can be seen, the ASVD has the smallest number of optimization calls followed by the ALS, the MASVD and the MALS.
One notices as well that the number of optimization calls for the two random data sets are close to each other for the respective algorithms. The real
data set takes most optimization calls, even though it probably profits from more potential correlations. However, the larger number of optimization calls may
also result from the different precision of one element of the third-order tensor (16bit vs. 8bit values).
Another explanation might be that it was difficult to find good rank one bases for a real data set (the error is approx. 70\% for the $512^3$ tensor). For
random data, the error stays around 63\%, probably due to a good distribution of the random values.
Otherwise, the number of optimization calls followed the same relationships as already seen in the timings measured for the rank one approximation algorithm.
For data sets larger than $128^3$, the time per optimization call stays roughly the same for any of the decomposition algorithms. However, the number of needed optimization calls is largest for the MALS and lowest for the ASVD.

\begin{figure}[htbp]
\begin{center}
 \includegraphics[width=1\columnwidth]{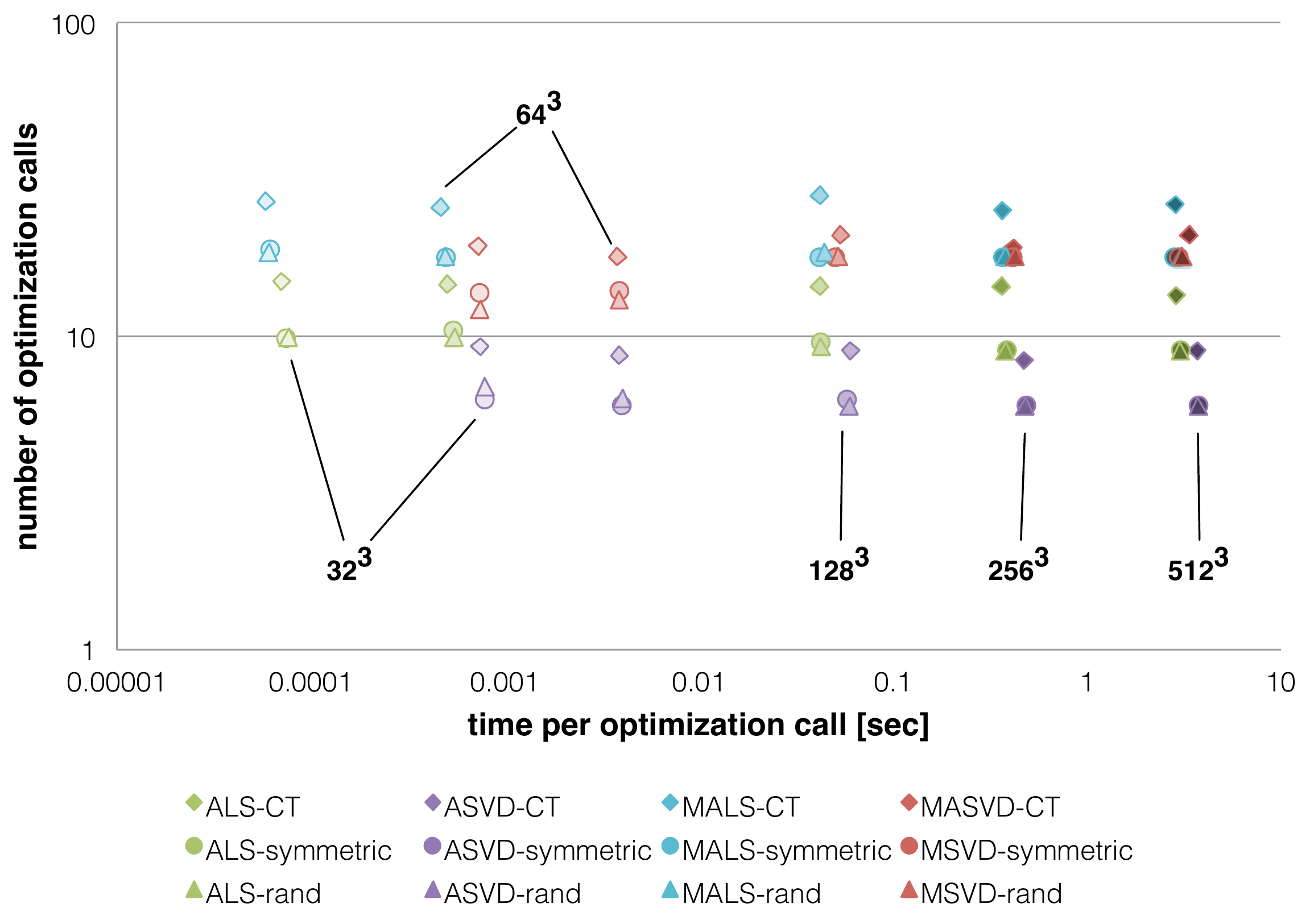}
\caption{Average time per optimization call put in relationship to the average number of optimization calls needed per algorithm and per data set  taken over 10 different initial random guesses.}\label{default1}
\end{center}
\label{fig:opt_calls}
\end{figure}

It is not only important to check how fast the different algorithms perform, but also what quality they achieve. This was measured by checking the
Frobenius norm of the resulting decompositions, which serves as a measure for the quality  of the approximation. In general, we can say that the higher
the Frobenius norm, the more likely it is that we find a global maximum. Accordingly, we compared the Frobenius norms in order to say whether the different
algorithms converged to the same stationary point.
In Figure 3, we show the absolute differences of the average Frobenius norms achieved by the ALS, ASVD, MALS and MASVD. The differences are taken with respect to the ALS.
As previously seen, the results for the real CT data set and the two random dataset differ. For the real data set, the differences for the achieved qualities are much smaller. Moreover, we see that the achieved quality for the ALS and the MALS are almost the same. A similar observation applies to the ASVD and the MASVD, which achieve almost the same quality.
We observed that all the algorithms reach the same stationary point for the smaller and medium data sets.
However, for the larger data sets ($ \ge 128^3$) the stationary points differ slightly. We suspect that either the same stationary point was not achieved, or the precision requirement of the convergence criterion was too high. That means that the algorithms stopped earlier,
since the results are not changing that much anymore
in the case that the precision tolerance for convergence is $0.0001$.

\begin{figure}[htbp]
\begin{center}
     \subfigure[CT data set]{\label{fig:sp_ct}
     \includegraphics[width=0.3\columnwidth]{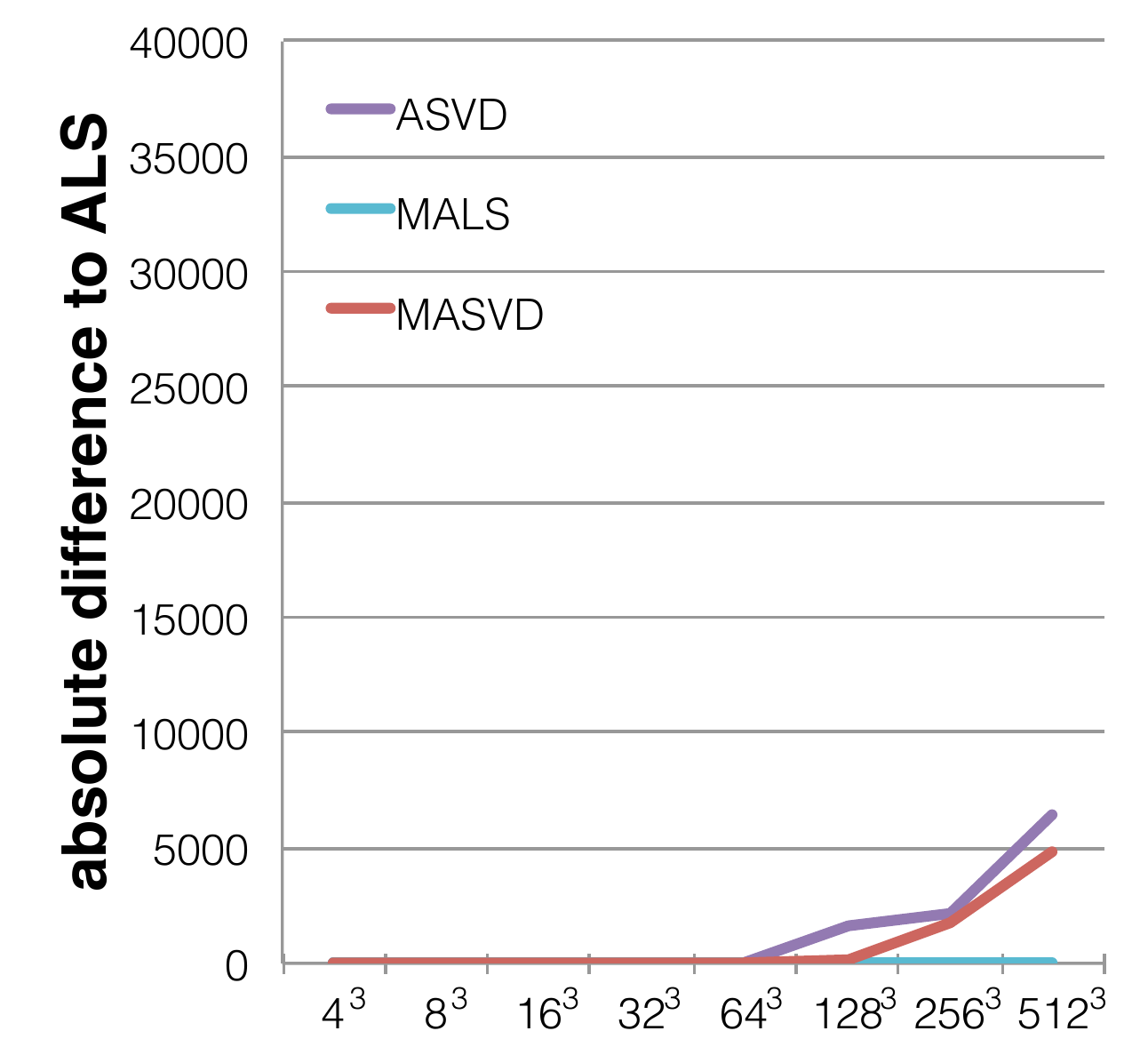}}
     \subfigure[symmetric data set]{\label{fig:sp_symm}
     \includegraphics[width=0.3\columnwidth]{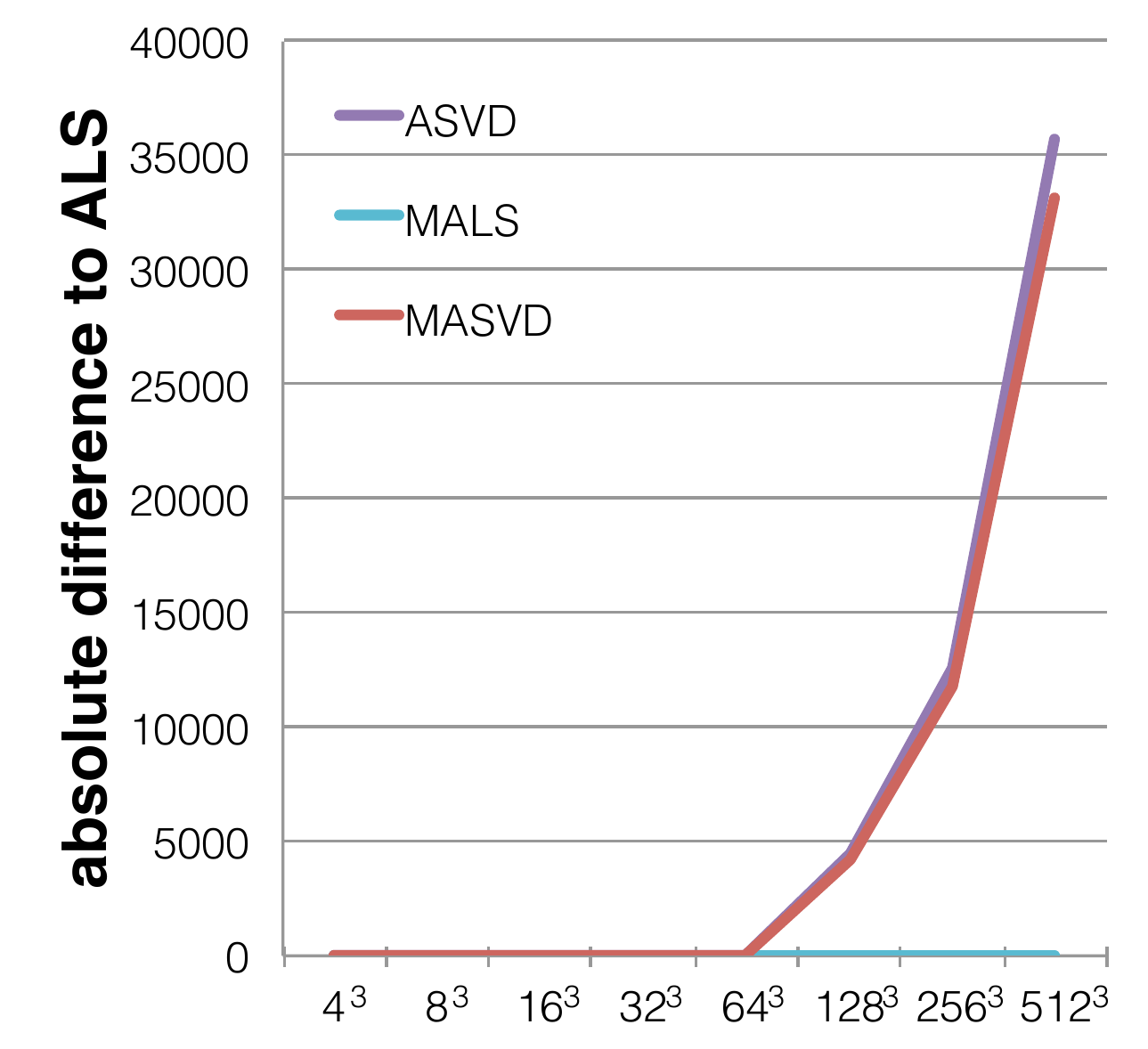}}
     \subfigure[random data set]{\label{fig:sp_rand}
     \includegraphics[width=0.3\columnwidth]{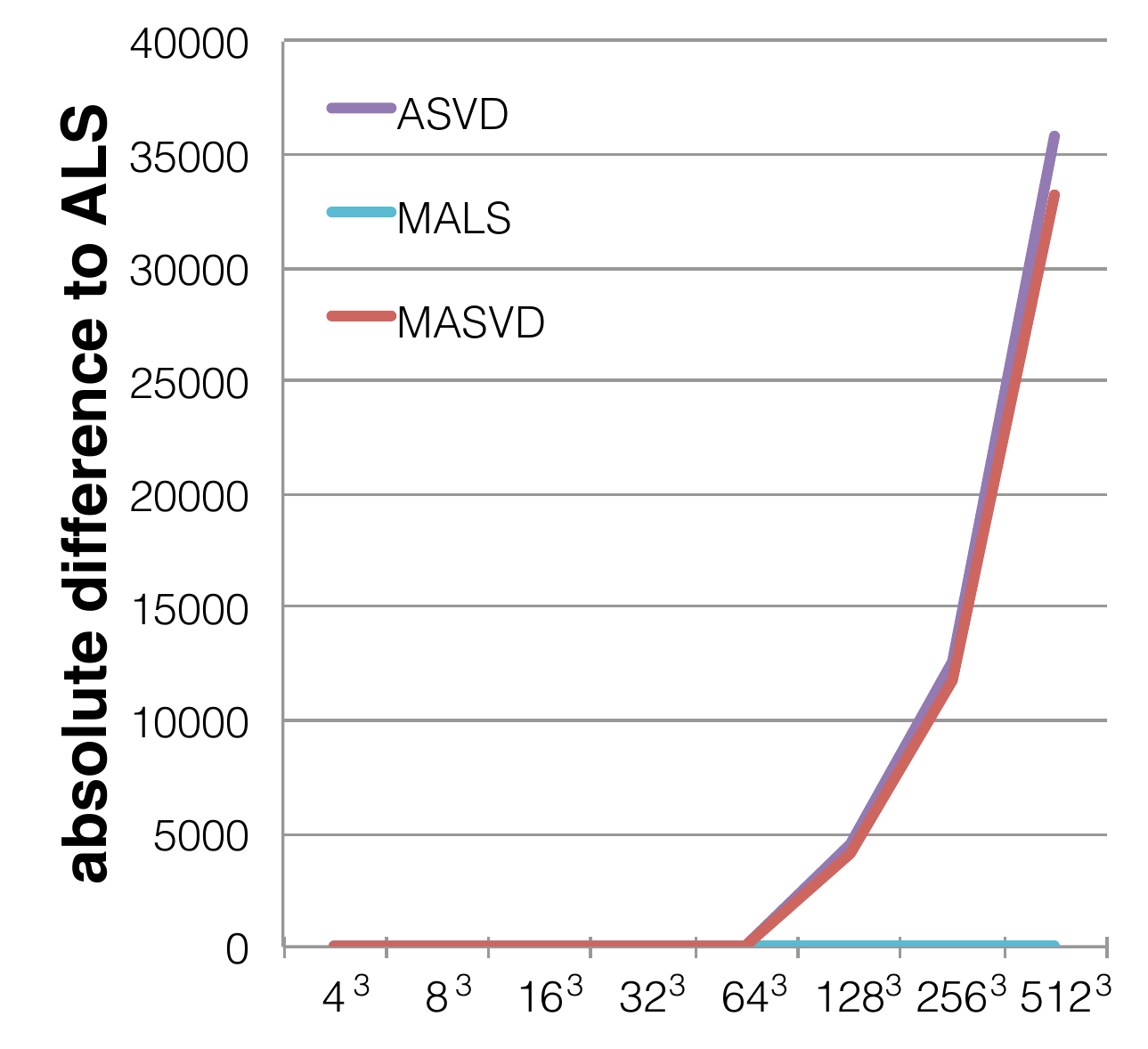}}
\caption{Differences of the achieved Frobenius norms by ALS, ASVD, MALS, and MASVD. The Frobenius norm of the approximations per algorithm and per data set are averages taken over 10 different initial random guesses.}
\label{fig:stat_point}
\end{center}
\end{figure}

Finally, the results of best rank one approximation for symmetric tensors using ALS, MALS, ASVD and MASVD show that
the best rank one approximation is also symmetric, i.e., is of the form $a\uu\otimes\vv\otimes \w$, where $\uu\approx \vv\approx \w \in \rS^{m-1}$.  This
confirms an observation made by Paul Van Dooren, (private communication), and the main result in \cite{Fri11}, which claims that the best rank one
approximation of a symmetric tensor can be always chosen symmetric.  The results of ASVD and MASVD give a better symmetric rank one approximation, i.e.,
$\uu-\vv,\uu-\w$ in ASVD and MASVD are smaller than in ALS and MALS.

\section{Conclusions}\label{sec:conc}

We have presented a new alternating algorithm for the computation of the best rank one approximation to a d-mode tensor. In contrast to the alternating
least squares method, this method uses a singular value decomposition in each step. In order to achieve guaranteed convergence to a semi-maximal point, we
have modified both algorithms.
%and also presented a  Newton type acceleration method.
We have run extensive numerical tests to show the performance and
convergence behavior of the new methods.

  \section*{Acknowledgements} The authors thank the OsiriX community for providing the MELANIX data set, and the referees for their comments.

 \bibliographystyle{plain}

\section*{Appendix: Remarks on local semi-maximality}
In this appendix we discuss the notion of an isolated critical point of a function $f$ which is semi-maximal but not maximal.
The main emphasize is to characterize semi-maximal points for which the alternating maximization iteration, abbreviated as AMI, converges
to the critical point at least for some nontrivial choices of the starting points.
We explain the convergence issues for ALS on local semi-maximality by the help of the AMI.

Consider a polynomial function $p(\bt),\ \bt\in\R^N$ and let $M\subset \R^N$ be a smooth compact manifold of dimension $L$.
Denote by $g$ the restriction of $p$ to $M$.  For example, in the three mode case we let $N=m+n+l$, $\bt=(\x,\y,\z)$, $p(\bt)=\cT\times
(\x\otimes\y\otimes\z)$ and  $M=\rS^{m-1}\times \rS^{n-1}\times \rS^{l-1},\ L=N-3$.
Assume that a point $\bt_{\star}\in M$ is a non-degenerate critical point
of $g$ on $M$.  We take local coordinates around $\bt_{\star}$, so
 that in these local coordinates $\bt_{\star}$ corresponds to
the zero vector of dimension $L$, denoted as $\0_L$.  So the open connected neighborhood of $\bt_\star$ is identified with an open connected neighborhood
$\0_L\in \R^L$.  Assume that the local coordinates around $\0_L$ are $\x\trans=[\x_1\trans,\ldots,\x_d\trans]\trans, \x_j\in \R^{m_j},j\in[d]$.

The AMI method consists of maximizing $g$ (or $f$) on $\x_j$ for $j=1,\ldots,d$, and then repeating the process. Let us discuss the details of the AMI for
a function $f$ given
by a quadratic form in the block  vector $\x\trans=[\x_1\trans,\ldots,\x_d\trans]\in \mathbb R^L$, given by
 \begin{eqnarray}\label{fformkv}
&& f=-\x\trans H\x,\
 H=\left[ \begin{array}{cccc} H_{1,1} &H_{1,2}& \cdots & H_{1,d}\\
 H_{2,1} & H_{2,2} & \ddots & H_{2,d} \\
 \vdots & \ddots & \ddots & \vdots\\
 H_{d,1} & \cdots & H_{d,d-1} & H_{d,d}
 \end{array} \right ],\\
 &&H_{p,q}\trans=H_{q,p},\ p,q\in[d]. \notag
 \end{eqnarray}
Note that locally we obtain this form for general $f$ via Taylor expansion and leaving off terms of order higher than two.

Consider the AMI iteration $\xi_{k-1}:=[\xi_{1,k-1}\trans,\ldots,\xi_{d,k-1}\trans]\trans\to \xi_{k}
=[\xi_{1,k}\trans,\ldots,\xi_{d,k}\trans]\trans\in\R^L$  for a function $f$
of the form (\ref{fformkv}) starting from a point $\xi_0$.  Then this iteration is the  block Gau{\ss}-Seidel iteration, see e.g. \cite{Var00}, applied to
the linear system $-H\xi=\xi_0$ with the block  symmetric matrix $H$, i.e.,
 \begin{equation}\label{GSit}
 -\sum_{\ell=1}^j H_{j,\ell}\xi_{\ell,k}=\sum_{\ell=j+1}^d H_{j,\ell}\xi_{\ell,k-1}, \quad j=1,\ldots,d, \quad k\in\N.
 \end{equation}
This iterative method can be expressed as $-L_H\xi_k=U_H\xi_{k-1}$, where $H=L_H+U_H$ is the decomposition of $H$ into the  block lower triangular part
$L_H$
and the strict  block upper triangular part $U_H$.  Assume that $L_H$ is invertible, which is equivalent to the requirement
that all diagonal blocks $H_{j,j}$ are invertible.
Then \eqref{GSit} is of the form $\xi_k=K\xi_{k-1}$, where
\begin{equation}\label{defK}
K:=-L_H^{-1}U_H.
\end{equation}
It is well known that an iteration $\xi_k=K\xi_{k-1}$ will converge to $\0_L$ for all starting vectors $\xi_0$ if and only if the spectral radius of $K$,
denoted as $\rho(K)$, is less than $1$.  If $\rho(K)\ge 1$ then the iteration will converge to $\0_L$ if and only if $\xi_0$ lies in the invariant subspace
of $K$ associated with the eigenvalues of modulus less than $1$.

Assume in the following that  $\0_L:=[\0_{m_1}\trans,\ldots,\0_{m_d}]\trans$
is a semi-maximal point, i.e., that
all diagonal blocks $H_{j,j},\ j\in [d]$ of $H$ are positive definite. Then it follows from a classical result of Ostrowski, see e.g. \cite[Thm
3.12]{Var00}, that the iteration \eqref{GSit} converges to $\0_L$ if and only if $H$ is positive definite, which is equivalent to $\rho(K)<1$.
Clearly, in this case $\0_L$ is non-maximal  for $f(\xi)$ if and only if  $H$ is indefinite.

We summarize these observations to give a precise condition on $\xi_0$ so that  the iteration~\eqref{GSit} converges to zero, which in the particular case discussed here can be proved easily. We give a proof for completeness.
\begin{theo}\label{AMIGS}
Let $\0_L:=[\0_{m_1}\trans,\ldots,\0_{m_d}]\trans$ be a semi-maximal point of $f(\xi)=-\xi^TH\xi$, i.e., each $H_{i,i}$ is positive definite and let $K$ be
given by \eqref{defK}.  Denote by $\alpha,\beta,\gamma$ the number of eigenvalues $\lambda$ of $K$, counting with multiplicities,
satisfying $|\lambda|<1, |\lambda|>1, |\lambda|=1$, respectively.
Assume that $H$ has $\pi, \nu,\zeta$ positive, negative and zero eigenvalues, respectively.
Then
\begin{eqnarray}\label{pilowbd}
&&\pi\ge \max\{m_j, j\in[d]\},\\
&&\alpha=\pi, \quad \beta=\nu, \quad \gamma=\zeta.\label{equlabc}
\end{eqnarray}
Furthermore, all $\gamma$ eigenvalues of $K$  on the unit circle correspond to a unique eigenvalue $1$
of geometric multiplicity $\gamma$.  The corresponding eigenvectors of $K$ are the eigenvectors of $H$ corresponding to the zero
eigenvalue.
 \end{theo}
\proof We first prove \eqref{pilowbd}.  Let $H_{i,i}$ be the diagonal block of maximal size $m_i$.  Let $\tilde H$ be a principal submatrix
of $H$ of order $m_i+1$ which has $H_{i,i}$ as its submatrix.  The Cauchy interlacing theorem \cite{HorJ85} implies that the eigenvalues of $\tilde H$
interlace with the eigenvalues of $H_{i,i}$.  Since all eigenvalues of $H_{i,i}$ are positive it follows that $\tilde H$ has at least $m_i$ positive
eigenvalues and
%Continuing this process of increasing by one the order of principle %submatrices of $H$ we deduce that $H$ has at least $m_i$.
hence, \eqref{pilowbd} holds.

To prove (\ref{equlabc}), assume first that $\zeta\ge 1$. But if  $\x$ is an eigenvector of $H$ corresponding to the eigenvalue $0$ then  $K\x=\x$.  Hence
$\gamma\ge \zeta$, and
$1$ is an eigenvalue of $K$ of geometric multiplicity at least $\zeta$.

Let $\V_0$ be the null space  of $H$.  Then $K$ restricted to $\V_0$ is the identity operator.
Consider the quotient space $\bQ:=\R^L/\V_0$.  Clearly, $K$ and $H$ induce linear operators $K_1,H_1:\bQ\to \bQ$, where $H_1$ is nonsingular with $\pi$
positive eigenvalues and $\nu$ negative eigenvalues.
Observe also that if $\y,\z\in \R^L$ and $\y-\z\in \V_0$ then $\y\trans H\y=\z\trans H\z$. Thus, it is enough to study the eigenvalues of $K_1$, which
corresponds to the case where $H$ is nonsingular, which we assume from now on.

Observe that the AMI does not decrease the value of $f(\xi)$.  Moreover, $f(\xi_k)=f(\xi_{k-1})$ if and only if $\xi_{k-1}=\0_L$.
Let us,  for simplicity of notation, consider the iteration $\xi_k=K\xi_{k-1}$ in the complex setting,
i.e., we consider $F(\xi)=-\xi^* H \xi$ ,where $\xi\in \C^L$.
All the arguments can also be carried out in the real setting, by considering pairs of complex conjugate eigenvalues and
the corresponding real invariant subspace associated with
the real and imaginary part of an eigenvector.

Let $\lambda$ be an eigenvalue of $K$ and let $\xi_0$ be the eigenvector to $ \lambda$. Then $F(\xi_1)=|\lambda|^2 F(\xi_0)> F(\xi_0)$ which implies
that $|\lambda| \neq 1$. (This implies that  the only eigenvalue of $K$
of modulus $1$ can be the eigenvalue $1$, which corresponds to the eigenvalue $0$ of $H$.)

Observe next, that if $H$ is positive definite, then $F(\xi_0)<0$ and the inequality $F(\xi_1)>F(\xi_0)$ yields that
$|\lambda|^2<1$, i.e., $\rho(K)<1$, which is Ostrowski's theorem.

>From now on we therefore assume that $H$ is indefinite and nonsingular.
Assume that $F(\xi_0)\ge 0$ and $\xi_0\ne \0_L$.  Then $F(\xi_k)$ is an increasing sequence which either diverges to
$+\infty$ or converges to a positive number.  Hence we cannot have  convergence $\xi_k\to\0_L$.
More precisely, we have convergence $\xi_k\to \0_L$ if and only if $F(\xi_k)\le 0$ for all $k\ge 0$.

Let $\U_0\subseteq\U_1\subset \C^L$ be the invariant subspaces of $K$ corresponding to the eigenvalues $0$ and the eigenvalues $\lambda$ of modulus less
than
$1$ of
$K$, respectively.  So $K\U_0\subset \U_0$ and $K|\U_0$ is nilpotent.  Let $l_0=\dim\U_0$.
We have that $F(\xi)\le 0$ for all $\xi\in\U$.  Let $\V_-,\V_+\subset\C^L$ be the eigen-subspaces corresponding to negative and positive eigenvalues of
$H$,  respectively.   So $\pi=\dim \V_+, \nu:=\dim \V_-$ and $\pi+\nu=L$.
 Consider $\W= \mbox{\rm Range}\ (K^L)$.  Then
\[
\U_0\cap \W= \{\0_L\},\; \dim \W=L-l_0, \; K\W=\W,\; \W+\U_0=\C^L.
\]
With $\W_+:=\W\cap \V_+$, then we have that $\dim \W_+\ge \pi -l_0$ and $K_1:=K|\W$ is invertible.
Setting $\W_j=K_1^{-j}W_+$, we have that  $\xi_j\in\W_j$, and $F(K^k\xi_j)\le 0$ for $k=0,\ldots,j$,  and clearly, $\dim \W_j=\dim\W_+$.
Since the space of $\dim \W_+$ subspaces in $\C^L$ is compact, there exists a subsequence of $\W_{j_k},k\in\N$ which converges to a $\dim\W_+$
 dimensional subspace $\X\subset\C^L$. This subspace corresponds to
 the invariant subspace of $K$ associated with eigenvalues satisfying $0<|\lambda|<1$,
 since $F(K^k\xi)\le 0$ for all $k\ge 0$ and $\xi\in\X$.  Thus, $\X\cap \U_0=\{\0_L\}$ and $\U_1=\X+\U_0$.  Note that $\dim\U_1=\dim\X+\dim\U_0\ge \pi$.
 Since $F(\xi)\le 0$ for each $\xi\in\U_1$ ,it follows that $\dim\U_1=\pi$, i.e., $\alpha=\pi$.

 As $\alpha+\beta=L$,  it then  follows that $\beta=L-\alpha=L-\pi=\nu$.
 \qed

 As an example, if we apply the  ALS method for finding the maximum of the trilinear form $\cT\times(\x\otimes\y\otimes \z)$ restricted to
 $(\x,\y,\z)\in M=\rS^{m-1}\times\rS^{n-1}\times\rS^{l-1}$, then this is just the AMI for the local quadratic form $g$.  It is well known that $g$ may
 have
 several critical  points, some of whom are strict local maxima and local semi-maxima see
  \cite[Example 2, p. 1331]{DeLDV00A}. The above analysis shows that the ALS may converge to each of these points for certain appropriate starting points.
 For a specific $\cT\in \R^{m\times n\times l}$ one can expect that the ALS iteration exhibits a complicated dynamics.
 Hence, it is quite possible that in some cases the ALS method
 will not converge to a unique critical point, see also \cite{DeLDV00A,KolB09,KroD80}.

 \end{document}